\documentclass[12pt,a4paper]{article}

\renewcommand{\theequation}{\arabic{section}.\arabic{equation}}

\newcommand{\bea}{\begin{eqnarray}}
\newcommand{\eea}{\end{eqnarray}}
\newcommand{\nn}{\nonumber}
\newcommand{\qedsymb}{{\em{Q.E.D.}}}

\newtheorem{thm}{Theorem}[section]
\newtheorem{dfn}[thm]{Definition}
\newtheorem{iden}[thm]{Identity}

\def\G{{\cal G}}
\def\Z{{\mathbf Z}}
\def\N{{\mathbf N}}
\def\R{{\mathbf R}}
\def\C{{\mathbf C}}
\def\pa{{\partial}}
\def\H{{\cal H}}

\begin{document}

\begin{titlepage} 
\vspace*{0.5cm}
\begin{center}
{\Large\bf  A note on a canonical  dynamical $r$-matrix}  
\end{center}

\vspace{1.5cm}

\begin{center}
{\large L. Feh\'er\footnote{Corresponding author, E-mail: lfeher@sol.cc.u-szeged.hu}
and B.G. Pusztai}
\end{center}
\bigskip
\begin{center}
Department of Theoretical Physics,  University of Szeged \\
Tisza Lajos krt 84-86, H-6720 Szeged, Hungary \\
\end{center}

\vspace{2.2cm}

\begin{abstract}  
It is well known that a classical dynamical $r$-matrix can be associated with
every finite-dimensional self-dual Lie algebra $\G$ by the
definition $R(\omega):= f(\mathrm{ad} \omega)$, where $\omega\in \G$ and
$f$ is the holomorphic function given by 
$f(z)=\frac{1}{2}\coth \frac{z}{2}-\frac{1}{z}$ for $z\in \C\setminus 2\pi i \Z^*$.
We present a  new, direct proof of the statement that
this canonical $r$-matrix satisfies the modified classical
dynamical Yang-Baxter equation on $\G$.
\end{abstract}

\end{titlepage}

\setcounter{section}{-1}
\section{Introduction}
\setcounter{equation}{0}

Dynamical generalizations of the Yang-Baxter equations
and the associated algebraic structures are in the focus 
of current interest due to their applications in 
the theory of integrable systems and other 
areas of mathematical physics and pure mathematics
(see \cite{ES1} for a review).
The present paper contains a detailed study of a particular dynamical 
$r$-matrix, which is an important special case 
of the classical dynamical
$r$-matrices introduced in \cite{EV}.

Let us recall that dynamical $r$-matrices in the sense of 
Etingof-Varchenko \cite{EV} are associated with any subalgebra $\H$
of any (complex or real) Lie algebra $\G$.
By definition, a dynamical $r$-matrix is a (holomorphic or smooth) 
$\G\otimes\G$-valued function on an open subset $\check\H^*$
of the dual space $\H^*$ of $\H$ subject to the following three 
conditions.
First, $r$ must satisfy  the modified classical dynamical Yang-Baxter equation
(mCDYBE): 
\begin{equation}
[r_{12}, r_{13}]+[r_{12}, r_{23}]+[r_{13}, r_{23}]+ 
T_j^1 \frac{\pa r_{23}}{\pa \omega_j} 
-T_j^2 \frac{\pa r_{13}}{\pa \omega_j}
+T_j^3 \frac{\pa r_{12}}{\pa \omega_j}
=\varphi,
\label{0.1}\end{equation}
where $\varphi$ is some constant, $\G$-invariant element of
$\G\wedge\G\wedge\G$.
The $\omega_j$ are coordinates on $\H^*$ with respect to a basis
$\{ T_j\}$ of $\H$, and the usual tensorial notations 
as well as the summation convention are used.
The second condition is that $(r+ r^T)$, where
$(X_a\otimes Y^a)^T= Y^a \otimes X_a$, is a $\G$-invariant constant.  
The third condition requires the map 
$r: \check \H^* \rightarrow \G\otimes \G$ to be equivariant 
with respect to the (coadjoint and adjoint) infinitesimal 
actions of $\H$ on the corresponding spaces. 
The mCDYBE  becomes the CDYBE for $\varphi =0$.

In most applications $\G$ is a simple Lie 
algebra and $\H$ is (a subalgebra of) a Cartan subalgebra. 
Another distinguished special case is is obtained by taking $\H:=\G$.
We consider this latter case, and allow $\G$ to be
any self-dual Lie algebra for which $\G^*$ can be identified
with $\G$ by means of an invariant scalar product 
$\langle\ ,\ \rangle$. 
We here study the dynamical $r$-matrix
given by the formula
\begin{equation}
r: \omega \mapsto r(\omega):= 
\langle T_j, f(\mathrm{ad} \omega ) T_k \rangle T^j \otimes T^k,
\qquad
\omega\in \check \G,
\label{0.2}\end{equation}
where $\check \G\subset \G$ is an open subset,
$\{T_j\}$ and $\{T^k\}$ denote dual bases of $\G$,
$\langle T_j, T^k\rangle = \delta_j^k$, 
and $f$ is the complex analytic function defined by 
\begin{equation}
f(z):=\frac{1}{2}\coth \frac{z}{2}- \frac{1}{z},
\qquad
z\in \C \setminus 2\pi i \Z^*.
\label{0.3}\end{equation}  
It is known that this $r$-matrix is a solution of the mCDYBE 
(\ref{0.1})  for $\H=\G\simeq \G^*$ with 
\begin{equation}
\varphi = -\frac{1}{4} f_{jk}^l T^j \otimes T^k \otimes T_l,
\qquad
[T_j,T_k]=f_{jk}^l T_l.
\label{0.4}\end{equation}
If $\G$ is a simple Lie algebra, then the mCDYBE for $r$ in (\ref{0.2})
follows from a general result (Theorem 3.14) in \cite{EV}.
Remarkably, this  $r$-matrix came to light 
naturally in two different applications, namely in the
context of equivariant-cohomology \cite{AM} and in the description 
of a Poisson structure on the chiral WZNW phase space 
compatible with classical $\G$-symmetry \cite{BFP}.  
A further reason for which
the $r$-matrix in (\ref{0.2}) is important is that it
can be reduced to certain self-dual subalgebras of $\G$, 
and thereby serves as a common `source' for a large family of dynamical  
$r$-matrices \cite{FGP1}.
We call it the canonical $r$-matrix of the self-dual Lie algebra $\G$.

The authors of \cite{AM} assumed $\G$ to be compact, while in \cite{BFP} 
$\G$ was taken to be a simple Lie algebra.  
In these papers the mCDYBE for the canonical
$r$-matrix was proved, independently, 
using the additional assumption that $\omega$ is near to zero,
so that $f(\mathrm{ad}\omega)$ is given with the aid of 
the power series expansion of $f(z)$ around $z=0$.
Though this is not obvious, the proofs 
found in \cite{AM,BFP} (see also \cite{AM+,Dub}) can in fact be adapted
to cover the case of a general self-dual Lie algebra as well.
In this case, a different proof of the mCDYBE appeared in \cite{ES2}. 
This proof is indirect and uses the restriction of $\omega$ 
to a neighbourhood of the origin.
The maximal domain of definition of $f(\mathrm{ad} \omega)$ 
contains all $\omega$ for which the eigenvalues
of $\mathrm{ad} \omega$ lie in $\C\setminus 2\pi i\Z^*$. 
Although the above-mentioned local proofs and the analyticity 
of $r(\omega)$ together imply the mCDYBE on this domain,
it could be enlightening to have an alternative direct proof, too.

The purpose of this paper is to present a direct proof
of the mCDYBE for the canonical $r$-matrix.
As opposed to the local power series expansion around $0$,
we here use the well known \cite{DS} holomorphic
functional calculus of linear operators to define 
$f(\mathrm{ad} \omega)$, and thus  
our proof is valid globally on the maximal domain of 
the `dynamical variable' $\omega$.
An advantage of our proof is that it also yields a
uniqueness result for the holomorphic function $f(z)$ that
enters the definition of the $r$-matrix in (\ref{0.2}).
Namely, by taking  formula (\ref{0.2}) 
as an ansatz the mCDYBE  translates into a functional
equation (eq.~(\ref{additional})) for the holomorphic function $f$ that 
has (\ref{0.3}) as its
unique solution under certain further natural conditions.
Despite its length, 
its elementary character and the uniqueness result that it implies
might justify presenting our proof. 
  
After this introduction, the paper consists of 2 sections and
3 appendices.
The proof of the mCDYBE is described in section 1.
It relies on some technical material collected in the appendices.
Appendix A is a recall of relevant basics of the functional calculus from
\cite{DS}, if necessary the reader might consult that first, the 
other two appendices should be consulted as referred to in the proof. 
Section 2 is devoted to a discussion of consequences 
of the proof, including the above-mentioned uniqueness result for 
the function $f$, and some comments.

\section{Proof of the mCDYBE for the canonical $r$-matrix}
\setcounter{equation}{0}

Let $\G$ be a finite-dimensional complex Lie algebra equipped with an 
invariant, symmetric, nondegenerate bilinear form $\langle\  ,\ \rangle$.
For the structure of such Lie algebras,  see e.g.~\cite{Figu}. 
We call these Lie algebras self-dual, 
since we identify $\G$ with $\G^*$ by means of the `scalar product' $\langle\ ,\ \rangle$.
Defining the transpose $A^T$ of an operator $A\in \mathrm{End}\left( \G \right)$ by 
$\langle A^{T} X,Y\rangle =\langle X,AY\rangle\ (\forall X,Y \in \G )$,
the invariance property of 
$\langle\ ,\ \rangle$ means that 
$\left( \mathrm{ad}\omega \right)^{T} 
=-\mathrm{ad}\omega \, ( \forall \omega \in \G )$, 
where $(\mathrm{ad}\omega)(X) =[ \omega ,X ]$.

Consider a map $r: \check \G \rightarrow \G\otimes \G$, where
$\check \G\subset \G$ is a nonempty open subset.
Then there exists a unique map $R: \check \G \rightarrow \mathrm{End}(\G)$ for which 
\begin{equation}
r(\omega)= \langle T_j, R( \omega ) T_k \rangle T^j \otimes T^k,
\qquad \forall \omega\in \check \G,
\end{equation}
where $\{T_j\}$ and $\{T^k\}$ denote dual bases of $\G$.
The directional derivatives of $R$ are given by
\begin{equation}
\left( \nabla_{S} R \right) \left( \omega \right) :=
\frac{\mathrm{d}}{\mathrm{d}t} \Bigg|_{t=0} 
R\left( \omega +tS\right) , \qquad
\forall S \in \G , \, \omega \in \check{\G}, 
\label{derdir}\end{equation}
and the `gradient' of $R$ is defined by
\begin{equation}
\left< X, \left( \nabla R\right) \left( \omega\right) Y\right> :=
 T^{j}\left< X, \left( \nabla_{T_{j}} 
R\right)\left(\omega\right) Y\right> ,
\qquad \forall X,Y\in\G,\, \omega\in\check{\G}.
\label{grad}\end{equation}

If $r$ is antisymmetric, i.e., $R^T(\omega)=-R(\omega)$,
then the mCDYBE (\ref{0.1}) for 
$r$ with $\varphi$ in (\ref{0.4})
is in fact equivalent to the following equation for $R$: 
\bea
\lefteqn{ \frac{1}{4} \left[ X,Y\right] 
+\left[ R\left(\omega\right) X,R\left(\omega\right)Y \right]
-R\left(\omega\right)
\left(
\left[ R\left(\omega\right) X,Y\right]
 +\left[ X,R\left(\omega\right) Y\right]
\right) } \nn \\
&&+\left< X,\left(\nabla R\right) \left(\omega\right) Y\right>
+\left( \nabla_{Y} R\right) \left(\omega\right) X
-\left( \nabla_{X} R\right) \left(\omega\right) Y=0,
\quad \forall X,Y\in\G,\,\omega\in\check{\G}  .
\label{cdybe}\eea
The $\G$-equivariance of the map $r: \check \G \rightarrow \G \otimes \G$
can be expressed as
\begin{equation}
(\nabla_{[S,\omega]} R)(\omega)= [\mathrm{ad} S, R(\omega)]
\qquad
\forall S\in \G,\, \omega\in \check \G.
\label{equivar}\end{equation}
 
After these remarks, we are ready to study the canonical $r$-matrix.  
From now on we use 
\begin{equation}
\check{ \G } :=
\left\{\, \omega \in  \G \mid 
\sigma \left(\mathrm{ad}\omega\right)\, \cap
2 \pi i \Z^{*} =\emptyset \, \right\}  ,
\end{equation}
which is a nonempty open subset in $\G$. 
Here and below $\sigma \left(\mathrm{ad}\omega\right)$
denotes the spectrum of 
$\mathrm{ad}\omega$ ($\forall \omega \in \G$),
and sometimes we use the notation 
$\bar{\omega} :=\mathrm{ad}\omega$ for brevity.
With the aid of the 
familiar  holomorphic functional calculus (see Appendix A), we can 
define an operator valued dynamical $r$-matrix  
$R\colon\check{\G} \to \mathrm{End}\left( \G \right)$ by 
\begin{equation}
\omega\mapsto R\left( \omega\right) :=f\left(\mathrm{ad}\omega\right) 
=\frac{1}{2\pi i} \int_C d\xi f(\xi) (\xi I - \mathrm{ad}\omega)^{-1},
\label{rmatrix}\end{equation}
where $f$ is given in (\ref{0.3}). 
The curve $C$ encircles each eigenvalue of $\mathrm{ad}\omega$ and
$I$ is the identity operator on $\G$.
Now our main theorem can be formulated as follows.

\medskip\noindent
{\bf Theorem 1.}
{\em The mapping (\ref{rmatrix}) with $f$ in (\ref{0.3}) 
defines an antisymmetric $r$-matrix which 
satisfies the equivariance  condition (\ref{equivar}) and 
the mCDYBE given by (\ref{cdybe}).}
\bigskip

The antisymmetry of the $r$-matrix follows from (\ref{rmatrix}) 
by using that $f$ is an odd function, and the equivariance condition
(\ref{equivar}) is also an immediate consequence of (\ref{rmatrix}) (cf.~(\ref{fc2})).
Before verifying (\ref{cdybe}),
we gather some useful information 
and lemmas that make the calculations easier.

Let $\omega$ be an arbitrary fixed element of $\check{\G}$. For every 
$\lambda\in\C$, let $b_{\lambda}:=\mathrm{ad}\omega -\lambda I=\bar{\omega} 
-\lambda I\in \mathrm{End}\left(\G\right)$.
Thanks to the derivation property of $\mathrm{ad}\omega$, 
the $b_{\lambda}$'s 
enjoy the identities
\begin{equation}
b_{\alpha +\beta}^{n} \left[ X,Y\right]=
\sum_{j=0}^{n}\
\left( \begin{array}{c} n\\ j\end{array}\right)
\left[ b_{\alpha}^{j} X,b_{\beta}^{n-j} Y\right],\quad
\forall X,Y\in\G,\,\forall\alpha ,\beta\in\C .
\label{propertyofb}\end{equation}
By means of the 
$\G =\oplus_{\lambda\in\sigma\left(\bar{\omega}\right)} N_{\lambda}$
Jordan decomposition, where 
$ N_{\lambda} =\mathrm{Ker}\left(b_{\lambda}^{\nu\left(\lambda\right)}\right)$
(see Appendix A), the $r$-matrix (\ref{rmatrix}) can be written as
\begin{equation}
R\left(\omega\right)=f\left(\bar{\omega}\right)=
\sum_{\lambda\in\sigma\left(\bar{\omega}\right)}
\sum_{k=0}^{\nu\left(\lambda\right) -1}
\frac{f^{\left( k\right)}\left(\lambda\right)}{k!}
b_{\lambda}^{k} E_{\lambda}.
\label{r_omega}
\end{equation}
We can regard this equation as the application of (\ref{fc3}) to the 
operator $\mathrm{ad}\omega$. 
Here $E_{\lambda}\in \mathrm{End}(\G)$ means the projection 
corresponding to the subspace $N_{\lambda}$.
Note also that $[N_\lambda, N_\mu]\subset N_{\lambda +\mu}$
is implied by (\ref{propertyofb}), with 
$N_\mu=\{ 0\}$ for any $\mu\notin \sigma(\bar \omega )$.

The mCDYBE (\ref{cdybe}) is linear in $X$ and $Y$. 
Therefore it is enough to 
prove this equation when $X\in N_{\lambda}$, $Y\in N_{\mu}$ are arbitrary 
elements of the subspaces associated with the eigenvalues 
$\lambda ,\mu\in\sigma\left(\bar{\omega}\right)$. So, from now on let 
$\lambda ,\mu$ be arbitrary, fixed eigenvalues of $\bar{\omega}$ and 
$X\in N_{\lambda}$, $Y\in N_{\mu}$ arbitrary, but fixed vectors. Applying 
the $r$-matrix (\ref{r_omega}) on these vectors, we obtain
\bea
R\left(\omega\right) X &=& f\left(\bar{\omega}\right) X=
\sum_{k=0}^{\nu(\lambda)-1} \frac{f^{\left( k\right)}\left(\lambda\right)}{k!}
b_{\lambda}^{k} X ,\nn\\
R\left(\omega\right) Y &=& f\left(\bar{\omega}\right) Y=
\sum_{l=0}^{\nu(\mu)-1} \frac{f^{\left( l\right)}\left(\mu\right)}{l!}
b_{\mu}^{l} Y . \label{rxry}
\eea

In the subsequent four lemmas we calculate the various terms of the mCDYBE (\ref{cdybe})
in a form that will prove convenient for verifying this equation.
In all expressions containing $\left[ b_{\lambda}^{k} X,b_{\mu}^{l} Y\right]$
it is understood that the indices $k, l$ vary as in (\ref{rxry}).

\medskip\noindent
{\bf Lemma 2.}
{\em 
If $\lambda ,\mu\in\sigma (\bar{\omega})$, $X\in N_{\lambda}$, $Y\in 
N_{\mu}$, then
\bea
\frac{1}{4}\left[X,Y\right]  &=&
\sum_{k,l} \lim_{(\alpha ,\beta )\to (\lambda ,\mu)} 
\frac{\pa^{k+l}}{\pa \alpha^{k} \pa \beta^{l}} \frac{1}{4}
\frac{\left[ b_{\lambda}^{k} X,b_{\mu}^{l} Y\right] }{k!l!} ,
\label{term1} \\
\left[ f\left( \bar{\omega} \right) X,f\left( \bar{\omega} \right) 
Y\right]  &=&
\sum_{k,l}\lim_{(\alpha ,\beta )\to (\lambda ,\mu )} 
\frac{\pa^{k+l}}{\pa \alpha^{k} \pa \beta^{l}}
f\left( \alpha\right) f\left( \beta\right)
\frac{\left[ b_{\lambda}^{k} X,b_{\mu}^{l} Y\right] }{k!l!} ,
\label{term2}\\
f\left( \bar{\omega} \right)
\left[ f\left( \bar{\omega}\right) X,Y\right]  &=&
\sum_{k,l}\lim_{(\alpha ,\beta ) \to (\lambda ,\mu )}
\frac{\pa^{k+l}}{\pa\alpha^{k} \pa\beta^{l}}
f\left( \alpha +\beta \right) f\left( \alpha \right)
\frac{\left[ b_{\lambda}^{k} X,b_{\mu}^{l} Y\right]}{k!l!} ,
\label{term3}\\
f\left( \bar{\omega}\right)
\left[ X,f\left( \bar{\omega}\right) Y\right]  &=&
\sum_{k,l}\lim_{(\alpha ,\beta )\to (\lambda ,\mu )}
\frac{\pa^{k+l}}{\pa\alpha^{k}\pa\beta^{l}}
f\left( \alpha +\beta\right) f\left(\beta\right)
\frac{\left[ b_{\lambda}^{k} X,b_{\mu}^{l} Y\right]}{k!l!} .
\label{term4}
\eea
}
\bigskip
\noindent 
{\em Proof.}  First, identity (\ref{d1}) from appendix C leads 
immediately  to (\ref{term1}) as 
\begin{equation}
\frac{1}{4} \left[ X,Y\right] =
\frac{1}{4} \left[ b_{\lambda}^{0} X,b_{\mu}^{0} Y\right]=
\sum_{k,l}\frac{\delta_{k,0}\delta_{l,0}}{4}
\frac{\left[ b_{\lambda}^{k} X,b_{\mu}^{l} Y\right]}{k!l!} =
\sum_{k,l}\frac{\pa^{k+l}}{\pa\alpha^{k}\pa\beta^{l}}
\frac{1}{4} \frac{\left[ b_{\lambda}^{k} X,b_{\mu}^{l} Y\right]}{k!l!}.
\end{equation}
Second, with the aid of (\ref{rxry}) and (\ref{d2}), we easily obtain (\ref{term2}) 
\bea
\left[ f\left(\bar{\omega}\right) X,f\left(\bar{\omega}\right) Y\right] 
&=&
\sum_{k,l} 
f^{\left( k\right)}\left( \lambda\right)
f^{\left( l\right)}\left( \mu\right) 
\frac{\left[ b_{\lambda}^{k} X,b_{\mu}^{l} Y\right]}{k!l!}\nn\\ 
&=&
\sum_{k,l}\lim_{(\alpha ,\beta )\to (\lambda ,\mu )}
\frac{\pa^{k+l}}{\pa\alpha^{k}\pa\beta^{l}}
f\left( \alpha\right) f\left( \beta\right) 
\frac{\left[ b_{\lambda}^{k} X,b_{\mu}^{l} Y\right]}{k!l!} .
\eea
Third, the calculation of 
\begin{equation}
f\left( \bar{\omega}\right)
\left[ f\left( \bar{\omega}\right) X,Y\right] =
f\left( \bar{\omega}\right)
\left[ \sum_{k}\frac{f^{\left( k\right)}\left(\lambda\right)}{k!}
b_{\lambda}^{k} X,Y\right] 
\end{equation}
goes as follows. 
Since $\left[\sum_{k} \frac{f^{\left( k\right)}\left(\lambda\right)}{k!} 
b_{\lambda}^{k} X,Y\right] \in N_{\lambda +\mu}$, (\ref{rxry})
yields
\bea
f\left(\bar{\omega}\right)
\left[ f\left(\bar{\omega}\right) X,Y\right] &=&
\sum_{k,l}
\frac{f^{\left( k\right)}\left(\lambda\right)
f^{\left( l\right)}\left(\lambda +\mu\right)}{k!l!}
b_{\lambda +\mu}^{l}
\left[ b_{\lambda}^{k} X, Y\right] \nn\\
&=&
\sum_{k,l}
\frac{f^{\left( k\right)}\left(\lambda\right)
f^{\left( l\right)}\left( \lambda +\mu\right)}{k!l!}
\sum_{j=0}^{l}
\left(\begin{array}{c} l\\j\end{array}\right)
\left[ b_{\lambda}^{k+l-j} X, b_{\mu}^{j} Y\right]  ,
\eea
where we used (\ref{propertyofb}). 
Introducing a new variable $s:=k+l$ for the 
summation, we have 
\bea
f\left( \bar{\omega}\right)
\left[ f\left( \bar{\omega}\right) X,Y\right] &=&
\sum_{s} \sum_{j=0}^{s} \sum_{l=j}^{s}
\left(\begin{array}{c} l\\ j\end{array}\right)
\frac{f^{\left( s-l\right)}\left(\lambda\right) 
f^{\left( l\right)}\left( \lambda +\mu\right)}{\left( s-l\right) !l!}
\left[ b_{\lambda}^{s-j} X, b_{\mu}^{j} Y\right] \\
&=&
\sum_{s} \sum_{j=0}^{s} \sum_{l=0}^{s-j}
\left( \begin{array}{c} l+j\\ j\end{array}\right)
\frac{f^{\left( j+l\right)}\left( \lambda +\mu\right)
f^{\left( s-j-l\right)}\left(\lambda\right)}
{\left( l+j\right) !\left( s-j-l\right) !}
\left[ b_{\lambda}^{s-j} X, b_{\mu}^{j} Y\right] \nn\\
&=&
\sum_{s} \sum_{j=0}^{s}
\frac{1}{j!\left( s-j\right) !}
\sum_{l=0}^{s-j}
\left(\begin{array}{c} s-j\\ l\end{array}\right)
f^{\left( j+l\right)}\left( \lambda +\mu\right)
f^{\left( s-j-l\right)}\left( \lambda\right)
\left[ b_{\lambda}^{s-j} X, b_{\mu}^{j} Y\right] .
\nn\eea
Using the Leibniz rule and introducing new 
summation variables as $l:=j$, $k:=s-j$, we obtain
\bea
f\left(\bar{\omega}\right)
\left[ f\left(\bar{\omega}\right) X,Y\right] &=&
\sum_{s} \sum_{j=0}^{s}
\frac{1}{j!\left( s-j\right) !}
\frac{\mathrm{d}^{s-j}}{\mathrm{d}\xi^{s-j}}\Bigg|_{\xi=\lambda}
f^{\left( j\right)}\left(\xi +\mu\right)
f\left(\xi\right) 
\left[ b_{\lambda}^{s-j} X, b_{\mu}^{j} Y\right] \nn\\
&=&
\sum_{k,l}\frac{\mathrm{d}^{k}}{\mathrm{d}\xi^{k}}\Bigg|_{\xi=\lambda}
f^{\left( l\right)}\left( \xi +\mu\right)
f\left(\xi\right) 
\frac{\left[ b_{\lambda}^{k} X, b_{\mu}^{l} Y\right]}{k!l!}.
\eea
By (\ref{d3}), this gives (\ref{term3}). 
Finally,  (\ref{term4}) is trivial consequence of (\ref{term3}). 
\qedsymb

\medskip\noindent
{\bf Lemma 3.}
{\em 
If $\lambda ,\mu\in\sigma\left(\bar{\omega}\right)$, $X\in N_{\lambda}$, 
$Y\in N_{\mu}$, then 
\begin{equation}
\left< X,\left( \nabla R\right) \left(\omega\right) Y\right> =
-\sum_{k,l}\lim_{(\alpha ,\beta )\to (\lambda ,\mu )}
\frac{\pa^{k+l}}{\pa\alpha^{k}\pa\beta^{l}}
\frac{f\left(\alpha\right) +f\left(\beta\right)}{\alpha +\beta} 
\frac{\left[ b_{\lambda}^{k} X, b_{\mu}^{l} Y\right]}{k!l!} .
\label{term5}\end{equation}
}
\bigskip
\noindent 
{\em Proof.} 
We obtain directly from the definitions (\ref{derdir}), (\ref{grad}),  
(\ref{rmatrix}) (see also (\ref{fc2})) that  
\begin{equation}
\left< X,\left(\nabla R\right)\left(\omega\right) Y\right> =
\frac{1}{ 2 \pi i} \int_{C} \mathrm{d} \xi f\left( \xi\right)
T^{j}\left< X,\rho_{\xi}\left(\bar{\omega}\right)
\left[ T_{j},\rho_{\xi}\left( \bar{\omega}\right) Y\right] \right>,
\end{equation}
where $\rho_\xi(\bar \omega)= (\xi I - \bar \omega)^{-1}$.
By using that 
$\rho_{\xi}\left(\bar{\omega}\right)^{T} = -\rho_{-\xi}\left(\bar{\omega}\right)$ 
and the invariance of 
$\left<\, ,\, \right>$, this expression is easily converted into 
\begin{equation}
\left< X,\left(\nabla R\right)\left(\omega\right) Y\right> 
=
\frac{1}{2\pi i}\int_{C}\mathrm{d}\xi f\left(\xi\right)
\left[ \rho_{-\xi}\left(\bar{\omega}\right) 
X,\rho_{\xi}\left(\bar{\omega}\right) Y\right] .
\end{equation}
We can apply the functional calculus to the holomorphic function 
$\rho_{\xi} : (\C\setminus \{\xi\})\rightarrow \C$ defined by 
$\rho_\xi: z \mapsto \left(\xi -z\right)^{-1}$. 
Thus we have 
\begin{equation}
\rho_{-\xi}\left(\bar{\omega}\right) X
=\sum_{k}\frac{\rho_{-\xi}^{\left( k\right)}\left(\lambda\right)}{k!}
b_{\lambda}^{k} X,\qquad
\rho_{\xi}\left(\bar{\omega}\right) Y
=\sum_{l} \frac{\rho_{\xi}^{\left( l\right)} \left( \mu\right)}{l!}
b_{\mu}^{l} Y,
\end{equation}
similarly to (\ref{rxry}).
Since $\rho_{-\xi}^{\left( k\right)}\left(\lambda\right) =k!\left(-\xi 
-\lambda\right)^{-\left( k+1\right)} =\left( -1\right)^{k+1} 
\rho_{\xi}\left( -\lambda\right)$, this leads to 
\begin{equation}
\left< X,\left(\nabla R\right)\left(\omega\right) Y\right> =
\sum_{k,l}\left(
\frac{ \left( -1\right)^{k+1}}{2\pi i}
\int_{C}\mathrm{d}\xi f\left( \xi\right)
\rho_{\xi}^{\left( k\right)}\left( -\lambda\right)
\rho_{\xi}^{\left( l\right)}\left(\mu\right) 
\right)
\frac{\left[ b_{\lambda}^{k} X, b_{\mu}^{l} Y\right]}{k!l!}  .
\label{integralformulaforterm5}\end{equation}
Now our task is to determine these integrals. Obviously, two 
different cases can appear. When $-\lambda =\mu$, the integrands have poles only at 
the point $\mu$. Alternatively, when $-\lambda\neq\mu$, the integrands have poles at 
the point $-\lambda$ and at the point $\mu$.

\noindent
\emph{The $\lambda +\mu =0$ case.}
In this case $\rho_{\xi}^{\left( k\right)}\left( -\lambda\right) 
\rho_{\xi}^{\left( l\right)}\left(\mu\right) =k!l!\left( \xi 
-\mu\right)^{-\left( k+l+1\right) -1}$.
Thanks to Cauchy's theorem, the 
integrals can be written as
\bea
\frac{\left( -1\right)^{k+1}}{2\pi i}
\int_{C}\mathrm{d}\xi f\left( \xi\right) 
\rho_{\xi}^{\left( k\right)}\left( -\lambda\right)
\rho_{\xi}^{\left( l\right)}\left(\mu\right) &=&
\frac{\left( -1\right)^{k+1} k!l!}{2\pi i}
\int_{C}\mathrm{d}\xi
\frac{f\left( \xi\right)}{ \left( \xi -\mu\right)^
{\left( k+l+1\right) +1}}= \nn\\ 
= \frac{\left( -1\right)^{k+1} k!l!}{\left( k+l+1\right) !}
f^{\left( k+l+1\right)}\left(\mu\right)
&=& -\lim_{(\alpha ,\beta )\to (\lambda ,\mu )}
\frac{\pa^{k+l}}{\pa\alpha^{k}\pa\beta^{l}}
\frac{f\left(\alpha\right) +f\left(\beta\right)}{\alpha +\beta} ,
\eea
where we used the identity (\ref{d5prime}). 
Thus (\ref{term5}) is valid in this case.

\noindent
\emph{The $\lambda +\mu\neq 0$ case.}
By $C_{\alpha}$ we denote a sufficiently small circle around the 
eigenvalue $\alpha\in\sigma\left(\bar{\omega}\right)$, which encircles 
this point in the positive sense. Using Cauchy's theorem in 
(\ref{integralformulaforterm5}), we can write
\bea
\lefteqn{ \frac{\left( -1\right)^{k+1}}{2\pi i}
\int_{C}\mathrm{d}\xi f\left(\xi\right)
\rho_{\xi}^{\left( k\right)}\left( -\lambda\right)
\rho_{\xi}^{\left( l\right)}\left(\mu\right) = }\nn\\
&&
=\left( -1\right)^{k+1}
\left\{
\frac{1}{2\pi i}\int_{C_{\mu}}\mathrm{d}\xi
f\left( \xi\right)
\rho_{\xi}^{\left( k\right)}\left( -\lambda\right)
\rho_{\xi}^{\left( l\right)}\left( \mu\right) +
\frac{1}{2\pi i}\int_{C_{-\lambda}}\mathrm{d}\xi
f\left(\xi\right)
\rho_{\xi}^{\left( l\right)}\left(\mu\right)
\rho_{\xi}^{\left( k\right)}\left( -\lambda\right)
\right\} \nn\\
&&
=\left( -1\right)^{k+1}
\left\{
\frac{\mathrm{d}^{l}}{\mathrm{d}\xi^{l}} \Bigg |_{\xi =\mu}
f\left(\xi\right)\left( -1\right)^{k+1}
\rho_{-\lambda}^{\left( k\right)}\left(\xi\right) +
\frac{\mathrm{d}^{k}}{\mathrm{d}\xi^{k}} \Bigg|_{\xi =-\lambda}
f\left(\xi\right)\left( -1\right)^{l+1}
\rho_{\mu}^{\left( l\right)}\left(\xi\right)
\right\} \nn\\
&&
=\left( -1\right)^{k}
\left\{
\sum_{a=0}^{l}\left( -1\right)^{k}
\left(\begin{array}{c} l\\ a\end{array}\right)
f^{\left( a\right)}\left(\mu\right)
\rho_{-\lambda}^{\left( k+l-a\right)}\left(\mu\right) +
\sum_{b=0}^{k}\left( -1\right)^{l}
\left(\begin{array}{c} k\\ b\end{array}\right)
f^{\left( b\right)}\left( -\lambda\right)
\rho_{\mu}^{\left( k+l-b\right)}\left( -\lambda\right)
\right\}\nn\\
&&
=-\left( -1\right)^{k+l}
\sum_{a=0}^{l}
\left(\begin{array}{c} l\\ a\end{array}\right)
\left( k+l-a\right) !\left( -1\right)^{a}
\frac{f^{\left( a\right)}\left(\mu\right)}{\left(\lambda 
+\mu\right)^{k+l+1-a}} \nn\\
&&\quad -\left( -1\right)^{k+l}
\sum_{b=0}^{k}
\left(\begin{array}{c} k\\ b\end{array}\right)
\left( k+l-b\right) !\left( -1\right)^{b}
\frac{ f^{\left( b\right)} \left(\lambda\right)}{ \left(\lambda 
+\mu\right)^{k+l+1-b}} .
\eea
Comparing this equation with (\ref{d5}),
we see  that when $\lambda +\mu\neq 0$ 
\begin{equation}
\frac{\left( -1\right)^{k+1}}{2\pi i}
\int_{C}\mathrm{d}\xi f\left(\xi\right)
\rho_{\xi}^{\left( k\right)}\left( -\lambda\right)
\rho_{\xi}^{\left( l\right)}\left(\mu\right) =
-\frac{\pa^{k+l}}{\pa\alpha^{k}\pa\beta^{l}}
\Bigg|_{(\alpha ,\beta )=(\lambda ,\mu )}
\frac{ f\left(\alpha\right) +f\left(\beta\right)}{\alpha +\beta} .
\end{equation}
Thus the proof of the lemma is complete. \qedsymb

\medskip\noindent
{\bf Lemma 4.}
{\em 
If $\lambda ,\mu\in\sigma\left(\bar{\omega}\right)$, $X\in 
N_{\lambda}$, $Y\in N_{\mu}$, then
\begin{equation}
\left(\nabla_{X} R\right)\left(\omega\right) Y=
\sum_{k,l}\lim_{(\alpha ,\beta )\to (\lambda ,\mu )}
\frac{\pa^{k+l}}{\pa\alpha^{k}\pa\beta^{l}}
\frac{f\left(\alpha +\beta\right) -f\left(\beta\right)}{\alpha}
\frac{\left[ b_{\lambda}^{k} X,b_{\mu}^{l} Y\right]}{k!l!}.
\label{term6}\end{equation}
}
\bigskip
\noindent 
{\em Proof.} 
As a consequence of (\ref{fc2}), the left hand side of 
(\ref{term6}) can be written as
\begin{equation}
\left(\nabla_{X} R\right)\left(\omega\right) Y=
\frac{1}{2\pi i}\int_{C}\mathrm{d}\xi f\left(\xi\right)
\rho_{\xi}\left(\bar{\omega}\right)
\left[ X,\rho_{\xi}\left(\bar{\omega}\right) Y\right].
\end{equation}
The application of the functional calculus 
(see also (\ref{term4}) and (\ref{d4})) gives
\begin{equation}
\rho_{\xi}\left(\bar{\omega}\right)
\left[ X,\rho_{\xi}\left(\bar{\omega}\right) Y\right] =
\sum_{k,l}\frac{\mathrm{d}^{l}}{\mathrm{d}\eta^{l}}
\Bigg|_{\eta =\mu}
\rho_{\xi}^{\left( k\right)}\left(\lambda +\eta\right)
\rho_{\xi}\left(\eta\right)
\frac{\left[ b_{\lambda}^{k} X,b_{\mu}^{l} Y\right]}{k!l!}.
\end{equation}
Therefore,
\begin{equation}
\left(\nabla_{X} R\right)\left(\omega\right) Y=
\sum_{k,l} \left\{
\sum_{j=0}^{l}\left(\begin{array}{c} l\\ j\end{array}\right)
\frac{1}{2\pi i}\int_{C}\mathrm{d}\xi f\left(\xi\right)
\rho_{\xi}^{\left( k+l-j\right)}\left(\lambda +\mu\right)
\rho_{\xi}^{\left( j\right)}\left(\mu\right) \right\}
\frac{\left[ b_{\lambda}^{k} X,b_{\mu}^{l} Y\right]}{k!l!} .
\label{integralformulaforterm6}\end{equation}
When $\lambda =0$, the integrands have poles only at the point $\mu$. 
If  $\lambda\neq0$, then the integrands have poles at the points 
$\lambda +\mu$ and $\mu$.

\noindent
\emph{The $\lambda =0$ case.}
In this case $\rho_{\xi}^{\left( k+l-j\right)}\left(\lambda +\mu\right) 
\rho_{\xi}^{\left( j\right)}\left(\mu\right) =\left( k+l-j\right) 
!j!\left(\xi -\mu\right)^{-\left(k+l+1\right)-1}$.
Thus  
\bea
\lefteqn{ \sum_{j=0}^{l}\left(\begin{array}{c} l\\ j\end{array}\right)
\frac{1}{2\pi i}\int_{C}\mathrm{d}\xi f\left(\xi\right) 
\rho_{\xi}^{\left( k+l-j\right)}\left(\lambda +\mu\right)
\rho_{\xi}^{\left( j\right)}\left(\mu\right) =} \nn\\
&&=
\sum_{j=0}^{l}\left(\begin{array}{c} l \\ j\end{array}\right)
\left( k+l-j\right) !j!
\frac{f^{\left( k+l+1\right)}\left(\mu\right)}{\left( k+l+1\right) !} =
\frac{k!l!f^{\left( k+l+1\right)}\left(\mu\right)}{\left( k+l+1\right) !}
\sum_{j=0}^{l}
\left( \begin{array}{c} \left( k+l\right) -j\\ 
\left( k+l\right) -l\end{array}\right) \nn\\
&&=
\frac{k!l!f^{\left( k+l+1\right)}\left(\mu\right)}{\left( k+l+1\right) !}
\left(\begin{array}{c} k+l+1\\ l\end{array}\right) =
\frac{f^{\left( k+l+1\right)}\left(\mu\right)}{k+1} \nn\\
&&=
\lim_{(\alpha ,\beta )\to (\lambda ,\mu )}
\frac{\pa^{k+l}}{\pa\alpha^{k}\pa\beta^{l}}
\frac{f\left(\alpha +\beta\right) -f\left(\beta\right)}{\alpha} ,
\eea
where we used the combinatorial identity (\ref{b2}) and (\ref{d6prime}). 
So in this case (\ref{term6}) holds. 

\noindent
\emph{The $\lambda\neq0$ case.}
Denote by $C_{\alpha}$ a sufficiently small circle around 
$\alpha \in \sigma\left( \bar{\omega}\right)$.
Then, by Cauchy's theorem, the relevant integrals in 
(\ref{integralformulaforterm6}) give
\bea
\lefteqn{
\frac{1}{2\pi i}\int_{C}\mathrm{d}\xi f\left(\xi\right)
\rho_{\xi}^{\left( k+l-j\right)}\left(\lambda +\mu\right)
\rho_{\xi}^{\left( j\right)}\left(\mu\right)=} \nn\\
&&=
\frac{1}{2\pi i}\int_{C_{\mu}}\mathrm{d}\xi f\left(\xi\right)
\rho_{\xi}^{\left( k+l-j\right)}\left(\lambda +\mu\right)
\rho_{\xi}^{\left( j\right)}\left(\mu\right) +
\frac{1}{2\pi i}\int_{C_{\lambda +\mu}}\mathrm{d}\xi f\left(\xi\right)
\rho_{\xi}^{\left( j\right)}\left(\mu\right)
\rho_{\xi}^{\left( k+l-j\right)}\left(\lambda +\mu\right) \nn\\
&&=
\frac{\mathrm{d}^{j}}{\mathrm{d}\xi^{j}}\Bigg|_{\xi =\mu}
f\left(\xi\right)
\rho_{\xi}^{\left( k+l-j\right)}\left(\lambda +\mu\right) +
\frac{\mathrm{d}^{k+l-j}}{\mathrm{d}\xi^{k+l-j}}
\Bigg|_{\xi =\lambda +\mu} f\left(\xi\right)
\rho_{\xi}^{\left( j\right)}\left(\mu\right) \nn\\
&&=
\left( -1\right)^{k+l-j+1}\sum_{a=0}^{j}
\left(\begin{array}{c} j\\ a\end{array}\right)
\left( k+l-a\right) !
\frac{f^{\left( a\right)} \left(\mu\right)}{\lambda^{k+l-a+1}} \nn\\
&&\quad +
\sum_{b=0}^{k+l-j}
\left(\begin{array}{c} k+l-j\\ b\end{array}\right)
\left( j+b\right) ! \left( -1\right)^{b}
\frac{f^{\left( k+l-j-b\right)}\left(\lambda +\mu\right)}{\lambda^{j+b+1}} .
\eea
Thus the coefficient of $\left[ b_{\lambda}^{k} X,b_{\mu}^{l} Y\right] 
/k!l!$ in (\ref{integralformulaforterm6}) is equal to the following 
expression:
\bea
\lefteqn{
\sum_{j=0}^{l}
\left(\begin{array}{c} l\\ j\end{array}\right)
\left\{ \left( -1\right)^{k+l-j+1}
\sum_{a=0}^{j}\left(\begin{array}{c} j\\ a\end{array}\right)
\left( k+l-a\right) !
\frac{f^{\left( a\right)}\left(\mu\right)}{\lambda^{k+l-a+1}}
\right. } \nn \\
&&+\left. \sum_{b=0}^{k+l-j}
\left(\begin{array}{c} k+l-j\\ b\end{array}\right)
\left( j+b\right) !\left( -1\right)^{b}
\frac{f^{\left( k+l-j-b\right)}\left(\lambda +\mu\right)}
{\lambda^{j+b+1}}\right\} .
\label{summation} \eea
Firstly, do the summation of the first part of (\ref{summation}):
\bea
\lefteqn{ \mathrm{Part1}(k,l):=
\sum_{j=0}^{l}\left(\begin{array}{c} l\\ j\end{array}\right)
\left( -1\right)^{k+l+1-j}
\sum_{a=0}^{j}\left(\begin{array}{c} j\\ a\end{array}\right)
\left( k+l-a\right) !
\frac{f^{\left( a\right)}\left(\mu\right)}{\lambda^{k+l-a+1}} =
} \nn \\
&&=
\left( -1\right)^{k+l+1}\sum_{a=0}^{l}
\frac{\left( k+l-a\right) !l!}{a!\left( l-a\right) !}
\frac{f^{\left( a\right)}\left(\mu\right)}{\lambda^{k+l-a+1}}
\left( -1\right)^{a}\sum_{j=0}^{l-a}
\left(\begin{array}{c} l-a\\ j\end{array}\right) \left( -1\right)^{j}
\nn \\
&&=
\left( -1\right)^{k+l+1}\sum_{a=0}^{l}
\frac{\left( k+l-a\right) !l!}{a!\left( l-a\right) !}
\frac{f^{\left( a\right)}\left(\mu\right)}{\lambda^{k+l-a+1}}
\left( -1\right)^{a}\delta_{l-a,0} 
=-\left( -1\right)^{k} k!
\frac{f^{\left( l\right)}\left(\mu\right)}{\lambda^{k+1}} .
\label{part1}\eea
Secondly, do the summation of the second part of (\ref{summation}). 
Introducing a new variable $m:=j+b$, we obtain
\bea
\lefteqn{\mathrm{Part2}(k,l):=
\sum_{j=0}^{l}
\left(\begin{array}{c} l\\ j\end{array}\right)
\sum_{b=0}^{k+l-j}
\left(\begin{array}{c} k+l-j\\ b\end{array}\right)
\left( j+b\right) !\left( -1\right)^{b}
\frac{f^{\left( k+l-j-b\right)}\left(\lambda +\mu\right)}{\lambda^{j+b+1}}
=} \nn\\
&&=
-\sum_{j=0}^{l}\sum_{b=0}^{k+l-j}\left( -1\right)^{m+1}
\frac{k!l!}{\left( k+l-m\right) !}
\frac{f^{\left( k+l-m\right)}\left(\lambda +\mu\right)}{\lambda^{m+1}}
\left( -1\right)^{j}
\left(\begin{array}{c} m\\ j\end{array}\right)
\left(\begin{array}{c} k+l-j\\ k\end{array}\right) \nn\\
&&=
-\sum_{m=0}^{l}\left( -1\right)^{m+1}
\frac{k!l!}{\left( k+l-m\right) !}
\frac{f^{\left( k+l-m\right)}\left(\lambda +\mu\right)}{\lambda^{m+1}}
\sum_{j=0}^{m}\left( -1\right)^{j}
\left(\begin{array}{c} m\\ j\end{array}\right)
\left(\begin{array}{c} k+l-j\\ k\end{array}\right) \nn\\
&&\,\,
-\sum_{m=l+1}^{k+l}\left( -1\right)^{m+1}
\frac{k!l!}{\left( k+l-m\right) !}
\frac{f^{\left( k+l-m\right)}\left(\lambda +\mu\right)}{\lambda^{m+1}}
\sum_{j=0}^{l}\left( -1\right)^{j}
\left(\begin{array}{c} m\\ j\end{array}\right)
\left(\begin{array}{c} k+l-j\\ k\end{array}\right)
\eea
By means of the combinatorial identities (\ref{b3}), (\ref{b4}),
we can simplify this formula. 
In fact, after a straightforward further computation,  we get
\begin{equation}
\mathrm{Part2}(k,l)=-\sum_{m=0}^{k}\left( -1\right)^{m+1}
\frac{k!}{\left( k-m\right) !}
\frac{f^{\left( k+l-m\right)}\left(\lambda +\mu\right)}{\lambda^{m+1}} .
\label{part2}\end{equation}
Now collecting equations 
(\ref{part2}), (\ref{part1}), (\ref{summation}), 
(\ref{integralformulaforterm6}), in the $\lambda\neq0$ case we can write
\bea
\left(\nabla_{X} R\right)\left(\omega\right) Y &=&
\sum_{k,l}\left\{\mathrm{Part1}(k,l)+\mathrm{Part2}(k,l)\right\}
\frac{\left[ b_{\lambda}^{k} X,b_{\mu}^{l} Y\right]}{k!l!} \nn\\
&=&
\sum_{k,l}\frac{\pa^{k+l}}{\pa\alpha^{k}\pa\beta^{l}}
\Bigg|_{(\alpha ,\beta )=(\lambda ,\mu )}
\frac{f\left(\alpha +\beta\right) -f\left(\beta\right)}{\alpha}
\frac{\left[ b_{\lambda}^{k} X,b_{\mu}^{l} Y\right]}{k!l!} ,
\eea
since equation (\ref{d6}) is valid. 
Hence lemma 4 is proved. \qedsymb

\medskip\noindent
{\bf Lemma 5.}
{\em 
If $\lambda ,\mu\in\sigma\left(\bar{\omega}\right)$, $X\in 
N_{\lambda}$, $Y\in N_{\mu}$, then
\begin{equation}
\left(\nabla_{Y} R\right)\left(\omega\right) X=
-\sum_{k,l}\lim_{(\alpha ,\beta )\to (\lambda ,\mu )}
\frac{\pa^{k+l}}{\pa\alpha^{k}\pa\beta^{l}}
\frac{f\left(\alpha +\beta\right) -f\left(\alpha\right)}{\beta}
\frac{\left[ b_{\lambda}^{k} X,b_{\mu}^{l} Y\right]}{k!l!} .
\label{term7}\end{equation}
}
\bigskip
\noindent 
{\em Proof.} 
This is a trivial consequence of the preceding lemma.

\medskip
Now we are in the position 
to verify the mCDYBE (\ref{cdybe}) 
for the canonical $r$-matrix (\ref{rmatrix}). 

\noindent
{\em Proof of theorem 1.} 
Let $\lambda ,\mu\in\sigma\left(\bar{\omega}\right)$ and
$X\in N_{\lambda}$, $Y\in N_{\mu}$. 
By applying the four lemmas, the 
left hand side of (\ref{cdybe}) can be written as
\bea
\lefteqn{
\frac{1}{4} \left[ X,Y\right] +
\left[ R\left(\omega\right) X,R\left(\omega\right) Y\right] -
R\left(\omega\right)\left(
\left[ R\left(\omega\right) X,Y\right] +
\left[ X,R\left(\omega\right) Y\right] 
\right)} \nn\\
&&\quad
 +\left< X,\left(\nabla R\right)\left(\omega\right) Y\right>
+\left(\nabla_{Y} R\right)\left(\omega\right) X
-\left(\nabla_{X} R\right)\left(\omega\right) Y= \nn\\
\lefteqn{ 
\quad =\sum_{k,l}\lim_{(\alpha ,\beta )\to (\lambda ,\mu )}
\frac{\pa^{k+l}}{\pa\alpha^{k}\pa\beta^{l}}
\left(
\frac{1}{4} +f\left(\alpha\right) f\left(\beta\right) -f\left(\alpha 
+\beta\right)\left( f\left(\alpha\right) +f\left(\beta\right)\right) 
\right.
} \nn\\
&&\quad\left.
-\frac{f\left(\alpha\right) +f\left(\beta\right)}{\alpha +\beta}
-\frac{f\left(\alpha +\beta\right) -f\left(\alpha\right)}{\beta}
-\frac{f\left(\alpha +\beta\right) -f\left(\beta\right)}{\alpha}
\right)
\frac{\left[ b_{\lambda}^{k} X,b_{\mu}^{l} Y\right]}{k!l!}. 
\label{proofid}\eea
This equals zero 
since the `addition formula' (\ref{additional}) is valid for the function
$f$ in (\ref{0.3}). \qedsymb

\section{Discussion}

We have shown that the canonical $r$-matrix defined by
(\ref{rmatrix}) with $f$ in (\ref{0.3}) satisfies the
mCDYBE (\ref{cdybe}). 
It is worth noticing that our proof implies a uniqueness result as well.
Suppose that we wish to define an antisymmetric 
solution of the mCDYBE (\ref{cdybe}) by the functional
calculus, i.e., by using {\em some} holomorphic complex function 
in formula (\ref{rmatrix}) now considered as an {\em ansatz}.  
For this formula to be well defined, the domain of holomorphicity 
of the function $f$ must contain zero, since this is always an
eigenvalue of $\mbox{ad} \omega$.
Moreover, for $R$ to be antisymmetric, which is in turn necessary for 
the equivalence of (\ref{cdybe}) to (\ref{0.1}) with $\varphi$ in (\ref{0.4}), 
$f$ must be an odd function.
Under these assumptions, the mCDYBE (\ref{cdybe}) for the ansatz  (\ref{rmatrix}) 
is in fact equivalent to the functional equation (\ref{additional}) 
for the unknown function $f$. 
Indeed, the whole calculation described in section 1 is valid 
for such an ansatz up to the equality in (\ref{proofid}).
The point then is that the functional equation (\ref{additional}) has a unique 
odd solution around the origin.
The proof of this statement is quite easy. 
By taking the $y\to0$ limit in (\ref{additional}) we obtain 
the differential equation for $f$ which appears in  (\ref{c1}). 
With the initial value $f\left( 0\right) =0$, which is implied by $f$ 
being odd, this differential equation has a unique, holomorphic solution around
the origin, namely the function 
$f\left( x\right) =\frac{1}{2}\coth\frac{x}{2}-\frac{1}{x}$.

So far we assumed the Lie algebra $\G$ to be complex, 
but the mCDYBE can be considered  for a real self-dual Lie algebra, too.
The real case arises naturally in applications \cite{AM,BFP}. 
Let us now suppose that $\G$ is the complexification of a real self-dual
Lie algebra, say $\G_r$.
Then it is not difficult to see that $R(\omega)$ given by (\ref{rmatrix})
maps $\G_r$ to $\G_r$ if $\omega\in \G_r$.
This is obviously the case if $\omega$ is near to zero,
where one can apply the power series expansion of $f$ around zero 
to define $R(\omega)$.
More generally, if $\omega\in \G_r$ then 
one may take the curve $C$ in (\ref{rmatrix}) to be invariant under complex 
conjugation as the eigenvalues of $\mathrm{ad}\omega$ occur in conjugate pairs.
By using this and $f(\bar z)= \bar f(z)$, 
complex conjugation of (\ref{rmatrix}) shows that 
$R(\omega) X \in \G_r$ if $\omega\in \G_r$ and $X\in \G_r$.
Thus the canonical $r$-matrix is a solution of 
the mCDYBE (\ref{cdybe}) in the real case as well.
 
Our use of the functional calculus,
which is  applicable 
to Banach algebras in general \cite{DS}, 
in the definition (\ref{rmatrix}) might serve as a starting point 
for future work towards  
generalizations of this canonical $r$-matrix  to certain 
infinite-dimensional self-dual Banach Lie algebras.
However, this represents a nontrivial problem since the 
above-presented proof of theorem 1
relies heavily on the finite-dimensionality of $\G$.

\bigskip
\bigskip
\noindent{\bf Acknowledgements.}
We are grateful to J.~Balog for reading the manuscript.
This investigation was supported in part by the Hungarian 
Scientific Research Fund (OTKA) under T034170 and M028418. 

\newpage

\renewcommand{\theequation}{\arabic{section}.\arabic{equation}}
\renewcommand{\thesection}{\Alph{section}}
\setcounter{section}{0} 

\section{Functional calculus of linear operators}
\setcounter{equation}{0}
\renewcommand{\theequation}{A.\arabic{equation}}

For convenient reference in the main text, 
in  this appendix we collect some result from the theory of 
bounded operators based  on chapter VII of the book \cite{DS}. 

 Let $X\neq\left\{ 0\right\}$ be a complex Banach space. 
The space of bounded linear operators on $X$ is 
denoted by $B\left( X\right)$, which is a Banach algebra in the usual way. 
Let $T\in B\left( X\right)$ be a bounded linear operator. 
The resolvent set of $T$ is given by 
${\cal R} \left( T \right) =\left\{ \, \lambda \in \C \mid \lambda I-T\, 
\mbox{invertible operator} \, \right\}$, where $I$ is the unit operator.
The spectrum $\sigma \left( T\right)$ of $T$ is the complement of 
${\cal R} \left( T\right)$. 
The formula ${\cal R}(T)\ni \xi \mapsto \rho_{\xi}\left( T\right) =\left( 
\xi I-T\right)^{-1}$ defines the resolvent function of $T$.  
Denote by $\mathcal{F} \left( T\right)$ the set of all complex functions $f$ 
that are holomorphic on some neighbourhood of $\sigma\left(T\right)$.
Then one can define the functions $f\left(T\right)$  of the operator $T$ as follows. 

\begin{dfn}
Let $f\in\mathcal{F}\left(T\right)$ and consider  a closed, rectifiable 
curve $C$ that lies in the domain of analiticity of $f$ and encircles the spectrum 
$\sigma\left(T\right)$ in the positive sense customary in the theory of 
complex variables. Then the operator $f\left(T\right)$ is defined by 
the equation
\begin{equation}
f\left(T\right)=\frac{1}{2\pi i}
\int\limits_{C}f\left(\xi\right)\rho_{\xi}\left(T\right)\mathrm{d}\xi.
\label{fc1}\end{equation}
\end{dfn}

It can be shown that $f\left(T\right)$ depends only on the function $f$, 
and not on the curve $C$.
Some important properties of this functional calculus are gathered in 
the following theorem.
\begin{thm}
If $f,g\in\mathcal{F}\left(T\right)$ and $\alpha ,\beta\in\C$ then
\begin{itemize}
\item $\alpha f+\beta g \in\mathcal{F}\left(T\right)$ and $\left(\alpha 
f+\beta g\right)\left(T\right)=\alpha f\left(T\right)
+\beta g\left(T\right)$,
\item $fg\in\mathcal{F}\left(T\right)$ and $\left(fg\right)\left(T\right)=
f\left(T\right) g\left(T\right)$,
\item if $f$ has the power series expansion 
$f\left(z\right)=\sum_{k=0}^{\infty}c_{k}z^{k}$ valid in a neighbourhood of 
$\sigma\left(T\right)$, then 
$f\left(T\right)=\sum_{k=0}^{\infty}c_{k}T^{k}$.
\end{itemize}
\end{thm}

One can define the directional derivatives,
$\left(\nabla_{S} f\right)\left(T\right)\in B(X)$, of $f\left(T\right)$ by 
\begin{equation}
\left(\nabla_{S} f\right)\left(T\right):=\frac{\mathrm{d}}{\mathrm{d}t}
\Bigg|_{t=0} f\left(T+tS\right), \qquad  S\in B\left(X\right).
\label{dirder}\end{equation}
The integral formula (\ref{fc1}) implies the equation 
\begin{equation}
\left(\nabla_{S} f\right)\left(T\right)=\frac{1}{2\pi i}\int\limits_{C} 
f\left(\xi\right)\rho_{\xi}\left(T\right) 
S\rho_{\xi}\left(T\right)\mathrm{d}\xi.
\label{fc2}\end{equation}

Now suppose that $X$ is a finite 
dimensional Banach space. 
In this case the spectrum $\sigma\left(T\right)$ of the operator $T$ has finitely many
elements, which are just the eigenvalues of $T$. 
The \emph{index} $\nu\left(\lambda\right)$ of an eigenvalue $\lambda$ 
is the smallest positive integer $\nu$ such that $\left(\lambda 
I-T\right)^{\nu}x=0$ for every vector $x$ for which $\left(\lambda 
I-T\right)^{\nu+1}x=0$. 
Introducing the invariant subspaces 
$N_{\lambda}:=\mathrm{Ker}\left(T-\lambda I\right)^{\nu\left(\lambda\right)} \, 
\left(\lambda\in\sigma\left(T\right)\right)$, one has the usual 
$X=\oplus_{\lambda\in\sigma\left(T\right)}N_{\lambda}$ Jordan 
decomposition of $X$.
\begin{thm}
If $dim\left(X\right)<\infty$ and $f\in\mathcal{F}\left(T\right)$, then
\begin{equation}
f\left(T\right)=\sum\limits_{\lambda\in\sigma\left(T\right)}
\sum\limits_{k=0}^{\nu\left(\lambda\right)-1}\frac{1}{k!}
f^{\left(k\right)}\left(\lambda\right)\left(T-\lambda I\right)^{k}
E_{\lambda},\label{fc3}\end{equation}
where $E_{\lambda}\in B(X)$ is the projection operator of the subspace 
$N_{\lambda}$.
\end{thm}

\section{Some combinatorial identities}
\setcounter{equation}{0}
\renewcommand{\theequation}{B.\arabic{equation}}

We here gather some elementary combinatorial identities needed in section 1.
\begin{iden}
If $k,l\in\N:=\{0,1,2,\ldots\}$, then
\begin{equation}
\sum\limits_{n=0}^{k}\left(-1\right)^{n} 
\frac{1}{n+l+1} 
\left( \begin{array}{c} k \\ n \end{array}\right)
=\frac{k!l!}{\left(k+l+1\right)!}.
\label{b1}\end{equation}
\end{iden}
\noindent 
{\em Proof.} 
By induction, with respect to k.
\begin{iden}
If $k,n\in\N$ and $0 \leq k\leq n$, then
\begin{equation}
\sum_{a=0}^{k}
\left( \begin{array}{c} n-a \\ n-k \end{array}\right)=
\left( \begin{array}{c} n+1 \\ k \end{array}\right).
\label{b2}\end{equation}\end{iden}
\noindent 
{\em Proof.} 
By induction with respect to n.
\begin{iden}
Let $k,l,m\in\N$ and $0\leq m\leq l$, then
\begin{equation}
\sum_{j=0}^{m} \left( -1\right)^{j}
\left( \begin{array}{c} m \\ j \end{array}\right)
\left( \begin{array}{c} k+l-j \\ k \end{array}\right)=
\left\{
\begin{array}{cl}
0 & \mbox{if $k<m$,}\\
\left(\begin{array}{c} k+l-m \\ l \end{array}\right) & \mbox{if $k\geq m$.}
\end{array}
\right.
\label{b3}\end{equation}
\end{iden}
\noindent 
{\em Proof.} 
Consider the smooth function
\begin{equation} 
\R\times\left(\R\setminus \left\{ 0\right\}\right) 
\ni\left(a,b\right)\mapsto
b^{k+l-m}\left(a+b\right)^{m} .
\end{equation}
Using the binomial theorem, we can write
\begin{equation}
b^{k+l-m} \left( a+b\right)^{m} =\sum_{j=0}^{m}
\left(\begin{array}{c} m\\ j \end{array}\right)
a^{j} b^{k+l-j}.
\label{ee1}\end{equation}
Let us differentiate this equation $k$-times with respect to  $b$. 
Then the left hand side gives 
\bea
\frac{\pa^{k}}{\pa b^{k}}\left( b^{k+l-m}\left(a+b\right)^{m}\right)&=&
\sum_{i=0}^{k}
\left(\begin{array}{c} k\\ i\end{array}\right)
\left( \frac{\pa^{k-i}}{\pa b^{k-i}}b^{k+l-m}\right)
\frac{\pa^{i}}{\pa b^{i}} \left(a+b\right)^m\nn\\
&=& \sum_{i=0}^{\mathrm{min}(m,k)}
\left(\begin{array}{c}k\\ i\end{array}\right)
\frac{\left(k+l-m\right)!\, m!}{\left( l+i-m\right) ! \left(m-i\right)!}
b^{l-m+i}\left(a+b\right)^{m-i}.
\eea
By evaluating this equation at $a=-1$, $b=1$, we obtain
\begin{equation}
\frac{\pa^{k}}{\pa b^{k}}\left( b^{k+l-m}\left(a+b\right)^{m}\right)
\Bigg|_{a=-1,\, b=1}=
\left\{ \begin{array}{cl}
0 &\mbox{if $k<m$,} \\
k!\left( \begin{array}{c} k+l-m \\ l\end{array}\right) &
\mbox{if $k\geq m$.} \end{array}\right.
\label{ee2}\end{equation}
At the same time, the right hand side of (\ref{ee1}) gives
\begin{equation} 
\frac{\pa^{k}}{\pa b^{k}}\sum_{j=0}^{m}
\left( \begin{array}{c} m\\ j\end{array}\right)
a^{j}b^{k+l-j}=k!\sum_{j=0}^{m}
\left(\begin{array}{c} m\\ j\end{array}\right)
\left(\begin{array}{c} k+l-j\\ k\end{array}\right)
a^{j}b^{l-j}. 
\end{equation}
It follows that
\begin{equation}
\frac{\pa^{k}}{\pa b^{k}}\sum_{j=0}^{m}
\left(\begin{array}{c} m\\ j\end{array}\right)
a^{j}b^{k+l-j}\Bigg|_{a=-1, \, b=1}=
k!\sum_{j=0}^{m}\left(-1\right)^{j}
\left(\begin{array}{c} m\\ j\end{array}\right)
\left(\begin{array}{c} k+l-j\\ k\end{array}\right) . 
\label{ee3}\end{equation}
Comparing (\ref{ee2}) and (\ref{ee3}) we see that our statement is valid. 
\qedsymb
\begin{iden}
Let $k,l,m\in\N$ and $l<m\leq k+l$, then
\begin{equation}
\sum_{j=0}^{l} \left(-1\right)^{j}
\left( \begin{array}{c} m\\ j\end{array} \right)
\left( \begin{array}{c} k+l-j\\ k\end{array} \right)=
\left\{ \begin{array}{cl}
0 &\mbox{if $k<m$,} \\
\left( \begin{array}{c} k+l-m \\ l \end{array} \right) &
\mbox{if $k\geq m$. }
\end{array}
\right.
\label{b4}\end{equation}
\end{iden}
\noindent 
{\em Proof.} 
Similar to the preceding identity. 

\section{Addition formula and further identities}
\setcounter{equation}{0}
\renewcommand{\theequation}{C.\arabic{equation}}

Let us consider the function 
$f\left(x\right)=\frac{1}{2}\coth\frac{x}{2}-\frac{1}{x}$.
This function is holomorphic on the whole complex plane 
except the points $2\pi i\Z^{*}$, where it has first order poles. 
Using the familiar 
$\coth x \coth y-\coth\left( x +y\right) 
\left( \coth x+\coth y\right)+1=0$ identity, the following   
`addition formula' can be obtained:
\begin{iden} If $x\neq 0$, $y\neq 0$, $x+y\neq 0$, then
\bea
\lefteqn{\frac{1}{4}+f\left( x\right) f\left( y\right) -
f\left( x+y\right) \left( f\left( x\right) +f\left( y\right) 
\right)}\nonumber\\
&&-\frac{f\left( x+y\right) -f\left( y\right)}{x}-
\frac{f\left( x+y\right) -f\left( x\right)}{y}-
\frac{f\left( x\right) +f\left( y\right)}{x+y}=0.
\label{additional}
\eea\end{iden}
On its domain of  holomorphicity, the function $f$ satisfies also 
the relations
\begin{equation}
f^{\left(k\right)}\left(-x\right)=\left(-1\right)^{k+1} 
f^{\left(k\right)}\left(x\right), \qquad
f'\left(x\right)+2\frac{f\left(x\right)}{x}+f^{2}\left(x\right)=
\frac{1}{4}.
\label{c1}\end{equation}
The first relation in (\ref{c1}) uses only the fact that $f$ is an odd function, 
while the second relation follows, for example, by taking the 
$y\to 0$ limit in (\ref{additional}). 
 
For convenience, 
we now collect some further identities that  
give the results for the differentiation of 
expressions of the type appearing in (\ref{additional}).
All these identities are obvious, and are actually valid    
for any odd holomorphic function $f$.
They are used in section 1 to derive 
the equality in (\ref{proofid}) for the $r$-matrix of the form 
in (\ref{rmatrix}). 

\begin{iden}
If $k,l\in\N=\{0,1,2,\ldots \}$, then
\bea
\frac{\pa^{k+l}}{\pa x^{k} \pa y^{l}} \frac{1}{4} &=& 
\frac{1}{4} \delta_{k,0} \delta_{l,0} , \label{d1} \\
\frac{\pa^{k+l}}{\pa x^{k} \pa y^{l}} 
f\left( x\right) f\left( y\right) &=& 
f^{\left(k\right)} \left( x\right)  f^{\left( l\right)} \left( y\right) ,
\label{d2} \\
\frac{\pa^{k+l}}{\pa x^{k} \pa y^{l}}
f\left( x+y\right) f\left( x\right) &=&
\frac{\mathrm{d}^{k}}{\mathrm{d} \xi^{k}} \Bigg|_{\xi=x}
f^{\left( l\right)}\left( \xi +y\right) f\left( \xi \right) ,
\label{d3} \\
\frac{\pa^{k+l}}{\pa x^{k}\pa y^{l}}
f\left( x+y\right) f\left( y\right) &=&
\frac{\mathrm{d}^{l}}{\mathrm{d}\xi^{l}}\Bigg|_{\xi =y}
f^{\left( k\right)}\left( \xi +x\right) f\left( \xi \right) .
\label{d4}
\eea \end{iden}
\begin{iden}
If $x+y\neq 0$, then
\bea
\frac{\pa^{k+l}}{\pa x^{k}\pa y^{l}}
\frac{f\left( x\right) +f\left( y\right)}{x+y} &&=
\left( -1\right)^{k+l}
\sum_{a=0}^{l}\left(\begin{array}{c} l\\ a\end{array}\right)
\left( k+l-a\right) ! \left( -1\right)^{a}
\frac{f^{\left( a\right)}\left( y\right)}{\left( x+y\right)^{k+l+1-a}}
\nonumber \\
&&\quad + \left( -1\right)^{k+l}
\sum_{b=0}^{k}\left(\begin{array}{c} k\\ b\end{array}\right)
\left( k+l-b\right) !\left( -1\right)^{b}
\frac{f^{\left( b\right)}\left( x\right)}{\left( x+y\right)^{k+l+1-b}}.
\label{d5}\eea
We  also have 
\begin{equation}
\lim_{x\to -y}
\frac{\pa ^{k+l}}{\pa x^{k}\pa y^{l}}
\frac{f\left( x\right) +f\left( y\right)}{x+y}=
\left( -1\right)^{k}
\frac{k!l!}{\left( k+l+1\right) !}
f^{\left( k+l+1\right) }\left( y\right) .
\label{d5prime}\end{equation}
\end{iden}
\noindent 
{\em Proof.} 
Equation (\ref{d5}) is a direct consequence of the Leibniz rule. 
To verify (\ref{d5prime}), let us introduce  $u:=x+y$,  $y=u-x$. 
By using power series expansion around $u=0$, we have  
\bea
\frac{f\left( x\right) +f\left( y\right)}{x+y} &&=
\frac{f\left( x\right) +f\left( u-x\right)}{u}=
\frac{f\left( x\right) -f\left( x-u\right)}{u}\nn\\
&&=\sum_{n=0}^{\infty} \left( -1\right)^{n}
\frac{f^{\left( n+1\right)} \left( x\right)}{\left( n+1\right) !}
u^{n}=
\sum_{n=0}^{\infty}
\left( -1\right)^{n}
\frac{f^{\left( n+1\right)}\left( x\right)}{\left( n+1\right) !}
\left( x+y\right)^{n} .
\eea
Differentiating this equation $l$-times with respect to $y$, we get that 
\begin{equation}
\frac{\pa^{l}}{\pa y^{l}}
\frac{f\left( x\right) +f\left( y\right)}{x+y}=\left( -1\right)^{l}
\sum_{n=0}^{\infty}\frac{\left( -1\right)^{n}}
{n!\left( n+l+1\right)}f^{\left( n+l+1\right)}\left( x\right)
\left( x+y\right)^{n} .
\end{equation}
Then differentiating $k$-times with respect to $x$, we obtain 
\bea
&&\frac{\pa^{k+l}}{\pa x^{k}\pa y^{l}}
\frac{f\left( x\right) +f\left( y\right)}{x+y} =
\left( -1\right)^{l}
\sum_{n=0}^{k}
\left(
\frac{\left( -1\right)^{n}}{n+l+1}
\sum_{j=0}^{n}
\left(\begin{array}{c} k\\ j\end{array}\right)
\frac{f^{\left( n+l+1+k-j\right)}\left( x\right)}{\left( n-j\right)!}
\left( x+y\right)^{n-j}\right) \nn\\
&&\qquad\qquad\quad+\left( -1\right)^{l}
\sum_{n=k+1}^{\infty} \left(
\frac{\left( -1\right)^{n}}{n+l+1}
\sum_{j=0}^{k}\left(\begin{array}{c} k\\ j\end{array}\right)
\frac{f^{\left( n+l+1+k-j\right)} \left( x\right)}{\left( n-j\right) !}
\left( x+y\right)^{n-j}
\right) .
\eea
Now, let us take the limit $x\to -y$. 
Using the combinatorial identity 
(\ref{b1}), we can see that 
\bea
\lim_{x\to -y}\frac{\pa^{k+l}}{\pa x^{k}\pa y^{l}}
\frac{f\left( x\right) +f\left( y\right)}{x+y} &&=
\left( -1\right)^{k}f^{\left( k+l+1\right)} \left( y\right)
\sum_{n=0}^{k}\frac{\left( -1\right)^{n}}{n+l+1}
\left(\begin{array}{c} k\\ n\end{array}\right)\nn\\
&&=\left( -1\right)^{k}
\frac{k!l!}{\left( k+l+1\right) !}
f^{\left( k+l+1\right)}\left( y\right),
\eea
whereby the proof is complete. \qedsymb
\begin{iden}
If $x\neq 0$, then
\begin{equation}
\frac{\pa^{k+l}}{\pa x^{k} \pa y^{l}}
\frac{f\left( x+y\right) -f\left( y\right)}{x} =
-\sum_{m=0}^{k}\frac{k!}{\left( k-m\right) !} \left( -1\right)^{m+1}
\frac{f^{\left( k+l-m\right)}\left( x+y\right)}{x^{m+1}} 
-\left( -1\right)^{k} k! 
\frac{ f^{\left( l\right)}\left( y\right)}{x^{k+1}} .
\label{d6}\end{equation}
In the limit case, we have
\begin{equation}
\lim_{x\to 0}
\frac{\pa^{k+l}}{\pa x^{k}\pa y^{l}}
\frac{f\left( x+y\right) -f\left( y\right)}{x}
=\frac{f^{\left( k+l+1\right)}\left( y\right)}{k+1}.
\label{d6prime}\end{equation}
\end{iden}
\noindent 
{\em Proof.} 
The verification of (\ref{d6}) is trivial. 
As for (\ref{d6prime}), 
the power series expansion of $f$ around $x=0$ implies that  
\begin{equation}
\frac{f\left( x+y\right) -f\left( y\right)}{x} =
\frac{1}{1!}f' \left( y\right) +\cdots 
+\frac{1}{\left( k+1\right) !} 
f^{\left( k+1\right)}\left( y\right) x^{k}
+\mathcal{O}\left( x^{k+1}\right) .
\end{equation}
By taking the derivatives of this equation, we obtain that  
\begin{equation}
\frac{\pa^{k}}{\pa x^{k}}
\frac{f\left( x+y\right) -f\left( y\right)}{x} =
\frac{f^{\left( k+1\right)} \left( y\right)}{k+1}
+\mathcal{O}\left( x\right), 
\end{equation}
and 
\begin{equation}
\frac{\pa ^{k+l}}{\pa x^{k}\pa y^{l}}
\frac{f\left( x+y\right) -f\left( y\right)}{x} =
\frac{f^{\left( k+l+1\right)}\left( y\right)}{k+1}
+\mathcal{O}\left( x\right),
\end{equation}
which implies (\ref{d6prime}).  \qedsymb
\begin{iden}
If $y\neq 0$, then
\begin{equation}
\frac{\pa^{k+l}}{\pa x^{k}\pa y^{l}}
\frac{f\left( x+y\right) -f\left( x\right)}{y} =
-\sum_{m=0}^{l} \frac{l!}{\left( l-m\right) !}
\left( -1\right)^{m+1} 
\frac{f^{\left( k+l-m\right)}\left( x+y\right)}{y^{m+1}} 
-\left( -1\right)^{l} l!
\frac{f^{\left( k\right)}\left( x\right)}{y^{l+1}}.
\label{d7}\end{equation}
In the limit case, 
\begin{equation}
\lim_{y\to0}
\frac{\pa^{k+l}}{\pa x^{k}\pa y^{l}}
\frac{f\left( x+y\right) -f\left( x\right)}{y} =
\frac{f^{\left( k+l+1\right)} \left( x\right)}{l+1} .
\label{d7prime}\end{equation}
\end{iden}
{\em Proof.}  This is an obvious consequence of the preceding identity.

\end{document}